\PassOptionsToPackage{frozencache}{minted}
\documentclass[11pt,a4paper]{amsart}
\usepackage{amsfonts, amsmath, amsthm, amssymb}
\usepackage[dvipsnames]{xcolor}
\usepackage{graphicx}
\usepackage{enumerate}
\usepackage{csquotes}
\usepackage{booktabs}
\usepackage[colorinlistoftodos, bordercolor=orange, backgroundcolor=orange!20, linecolor=orange, textsize=scriptsize]{todonotes}

\graphicspath{ {images/} }

\newtheorem{theorem}{Theorem}
\newtheorem{proposition}[theorem]{Proposition}
\theoremstyle{remark}

\usepackage[style=numeric,
   giveninits=true,
   sorting=nty,
   backend=biber,
   maxbibnames=3,
   isbn=false,
   backref=false]{biblatex}
\DeclareSourcemap{
  \maps[datatype=bibtex]{
    \map{
      \step[fieldsource=doi, final]
      \step[fieldset=url, null]
    }
  }
}

\newcommand{\CC}{\mathbb{C}}
\newcommand{\PP}{\mathbb{P}}
\newcommand{\QQ}{\mathbb{Q}}
\newcommand{\RR}{\mathbb{R}}

\newcommand{\field}{\Bbbk}

\usepackage{xspace}
\newcommand{\oscar}{\texttt{OSCAR}\xspace}
\newcommand{\exacus}{\texttt{EXACUS}\xspace}
\newcommand{\hecke}{\texttt{Hecke.jl}\xspace}
\newcommand{\flint}{\texttt{FLINT}\xspace}
\newcommand{\julia}{\texttt{Julia}\xspace}
\newcommand{\msolve}{\texttt{msolve}\xspace}
\newcommand{\TikZ}{\texttt{TikZ}\xspace}
\newcommand{\homcon}{\texttt{HomotopyContinuation.jl}\xspace}

\newcommand{\mathematica}{\texttt{Mathematica}\xspace}
\newcommand{\cgal}{\texttt{CGAL}\xspace}

\newcommand{\qepcad}{\texttt{QEPCADB}\xspace}
\newcommand{\macaulay}{\texttt{Macaulay2}\xspace}
\newcommand{\algsolving}{\texttt{Alge\-braic\-Solving.jl}\xspace}
\newcommand{\curvetopology}{\texttt{CurveTopology.jl}\xspace}
\newcommand{\pace}{\texttt{PACE.jl}\xspace}
\newcommand{\calcium}{\texttt{Calcium}\xspace}
\newcommand{\plotsjl}{\texttt{Plots.jl}\xspace}
\newcommand{\mcurve}{$M$-curve\xspace}
\newcommand{\mcurves}{$M$-curves\xspace}
\newcommand\scheme[1]{\langle #1 \rangle}
\usepackage[draft=true]{minted2}
\DeclareMathOperator{\conv}{conv}
\setminted{breaklines=true, breakanywhere=true}
\definecolor{lbcolor}{rgb}{0.9,0.9,0.9}
\setminted{bgcolor=lbcolor}
\setminted{fontsize=\footnotesize}

\usepackage[colorlinks=true, pdfstartview=FitV, linkcolor=Blue, citecolor=OliveGreen, urlcolor=BrickRed, breaklinks=false]{hyperref}

\title{A Known $M$-Curve of Degree $8$ in \oscar}
\hypersetup{pdftitle={A Known M-Curve of Degree 8 in OSCAR}, pdfauthor={Anne Frühbis-Krüger, Michael Joswig, Lars Kastner}}
\author[Frühbis-Krüger, Joswig and Kastner]{Anne Frühbis-Krüger \and Michael Joswig \and Lars Kastner}

\address[A. Frühbis-Krüger]{Universität Oldenburg, Institute of Mathematics}
\email{anne.fruehbis-krueger@uni-oldenburg.de}

\address[M.~Joswig, L.~Kastner]{Technische Universität Berlin, Chair of Discrete Mathe\-ma\-tics/Geo\-me\-try}
\email{\{joswig,kastner\}@math.tu-berlin.de}

\begin{document}

\begin{abstract}
  We show how the computer algebra system \oscar can be used to retrace geometric constructions of plane algebraic curves and provide topologically correct and visually pleasing drawings.
  The approach is illustrated on Shustin's construction from 1987 of a particular \mcurve of degree eight.
  We arrive at an explicit homogeneous trivariate polynomial with rational coefficients, and we provide a computational proof of all claimed properties of Shustin's original article.
\end{abstract}

\maketitle

\section{Introduction}

Hilbert's 16th problem asks about a topological classification of the real plane (projective) algebraic curves by degree.
Harnack~\cite{1876:Harnack:VieltheiligkeitEbenen} proved that the number of connected components (in the Euclidean topology) of a real curve of degree $d$ is bounded by $ M = \tfrac{1}{2}(d-1)(d-2)+1 $.
Curves which attain that bound are known as \emph{\mcurves}.
A large body of work from the 1980s culminated in a complete list of all isotopy types up to and including degree seven.
The reader is referred to articles by Ilyashenko \cite{Hilbert16} and Viro \cite{Viro:2008} for an overview.
The first case where the classification is not fully complete is degree $d=8$; see Orevkov \cite{2002:Orevkov:CurvesIsotopy} for details.

\begin{figure}[p]
   \includegraphics[width=.49\textwidth]{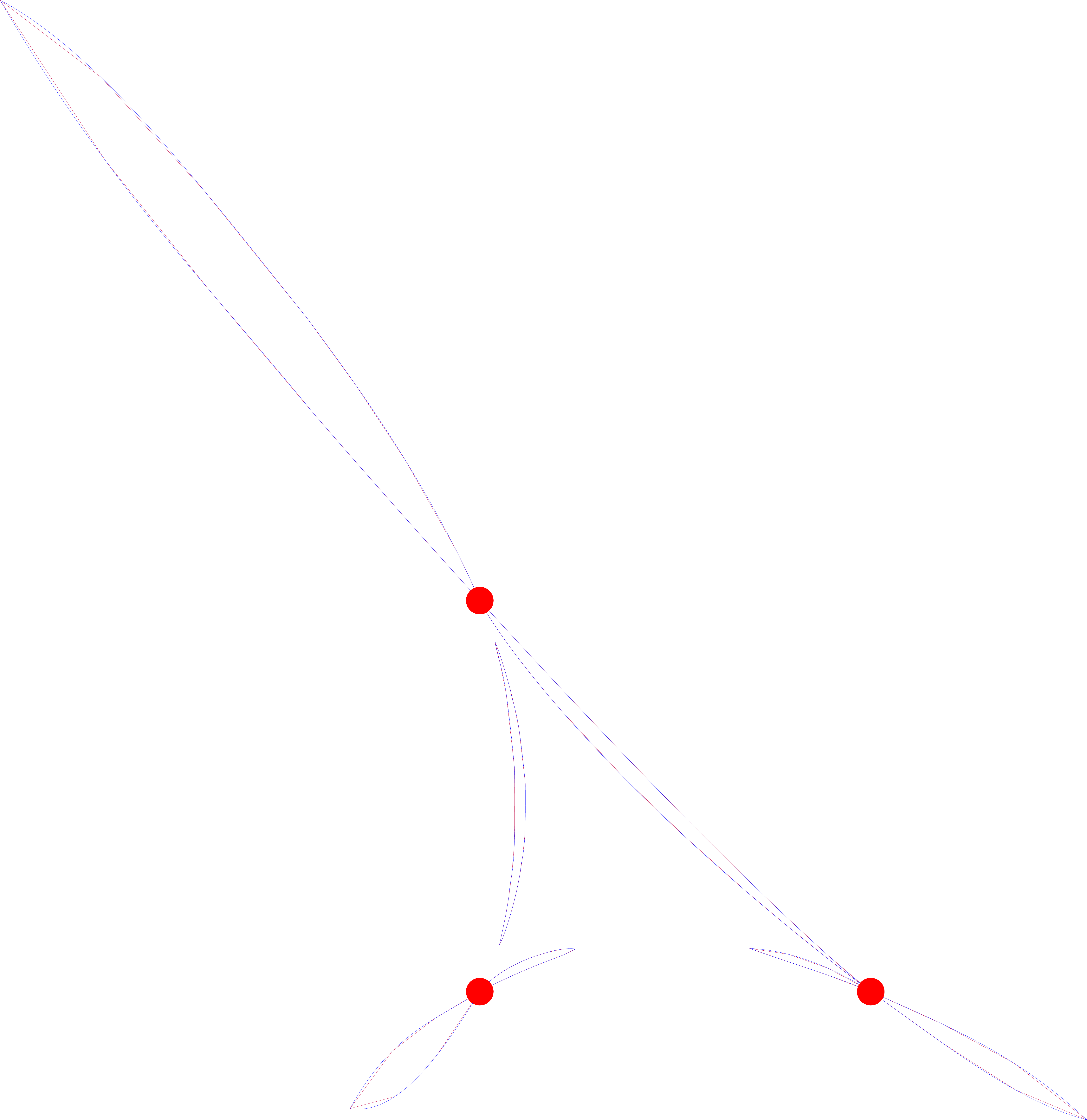}
   \hfill
   \includegraphics[width=.49\textwidth]{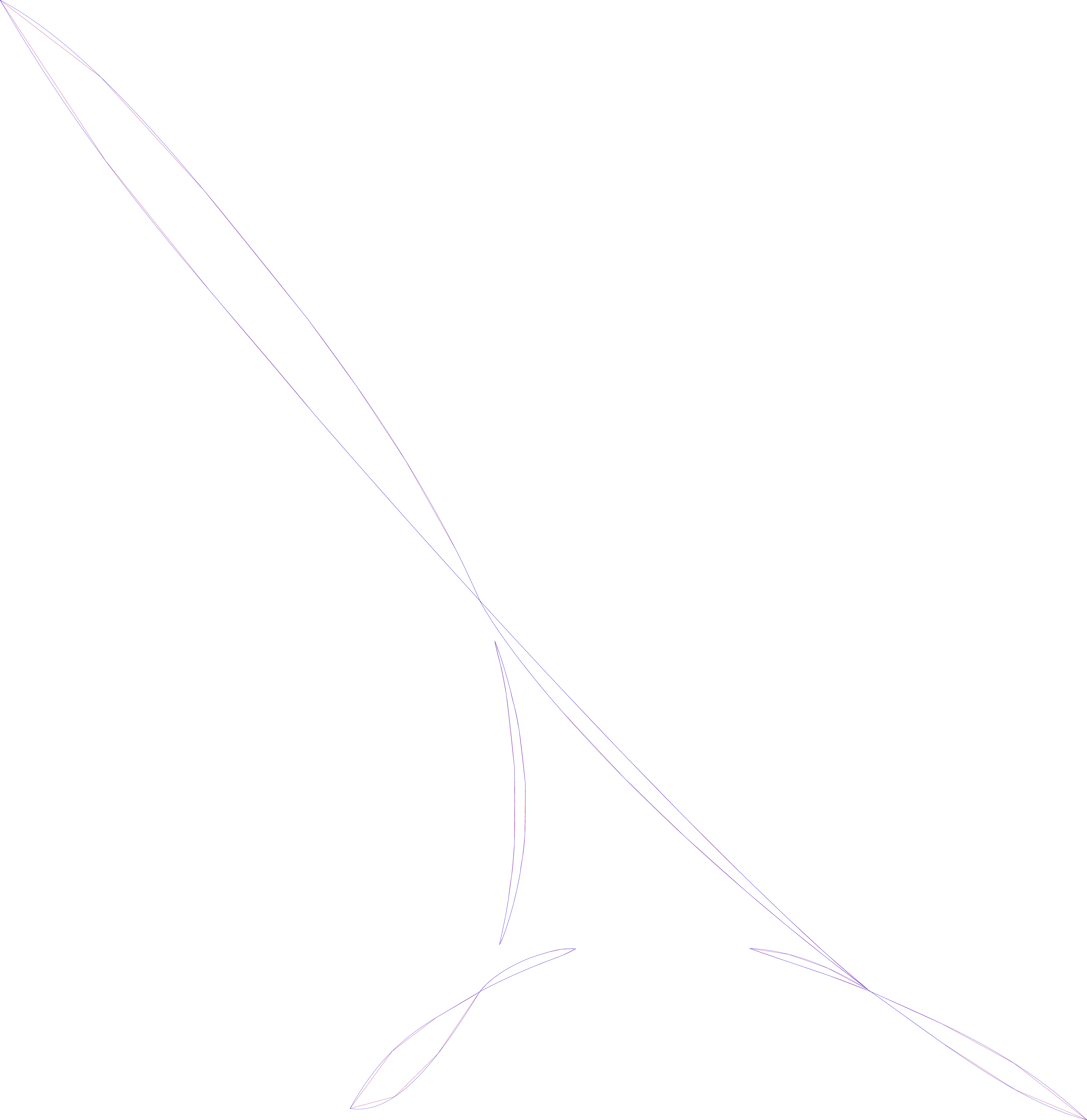}
   \caption{The real drawings of the singular curve $Q(C_4)$ (left) and the smooth curve $C_8$ (right) barely differ to the naked eye, and the topological change is invisible.}
   \label{fig:shustin_real}
   
   \vspace{1.4cm}
   \begin{tikzpicture}[scale=.6]
      \draw[black] (1.0, 7.5) to[out=270, in=180] (1.333333, 7.0);
\draw[black] (1.666667, 7.0) -- (1.333333, 7.0);
\draw[black] (2.0, 7.0) -- (1.666667, 7.0);
\draw[black] (1.0, 7.5) to[out=90, in=180] (1.333333, 8.0);
\draw[black] (1.666667, 8.0) -- (1.333333, 8.0);
\draw[black] (2.0, 8.0) -- (1.666667, 8.0);
\draw[black] (2.333333, 8.0) -- (2.0, 8.0);
\draw[black] (2.666667, 8.0) -- (2.333333, 8.0);
\draw[black] (3.0, 8.0) -- (2.666667, 8.0);
\draw[black] (2.333333, 7.0) -- (2.0, 7.0);
\draw[black] (2.666667, 7.0) -- (2.333333, 7.0);
\draw[black] (3.0, 7.0) -- (2.666667, 7.0);
\draw[black] (2.0, 1.5) to[out=270, in=180] (2.333333, 1.0);
\draw[black] (3.0, 1.5) to[out=180, in=0] (2.666667, 1.0);
\draw[black] (2.666667, 1.0) -- (2.333333, 1.0);
\draw[black] (2.0, 1.5) to[out=90, in=180] (2.333333, 2.0);
\draw[black] (3.0, 1.5) to[out=180, in=0] (2.666667, 2.0);
\draw[black] (2.666667, 2.0) -- (2.333333, 2.0);
\draw[black] (3.333333, 8.0) -- (3.0, 8.0);
\draw[black] (4.0, 7.5) to[out=180, in=0] (3.666667, 8.0);
\draw[black] (3.666667, 8.0) -- (3.333333, 8.0);
\draw[black] (3.333333, 7.0) -- (3.0, 7.0);
\draw[black] (4.0, 7.5) to[out=180, in=0] (3.666667, 7.0);
\draw[black] (3.666667, 7.0) -- (3.333333, 7.0);
\draw[black] (3.0, 1.5) to[out=0, in=180] (3.333333, 1.0);
\draw[black] (3.666667, 1.0) -- (3.333333, 1.0);
\draw[black] (4.0, 1.0) -- (3.666667, 1.0);
\draw[black] (3.0, 1.5) to[out=0, in=180] (3.333333, 2.0);
\draw[black] (3.666667, 2.0) -- (3.333333, 2.0);
\draw[black] (4.0, 2.0) -- (3.666667, 2.0);
\draw[black] (4.333333, 1.0) -- (4.0, 1.0);
\draw[black] (4.666667, 1.0) -- (4.333333, 1.0);
\draw[black] (5.0, 1.0) -- (4.666667, 1.0);
\draw[black] (4.333333, 2.0) -- (4.0, 2.0);
\draw[black] (4.666667, 2.0) -- (4.333333, 2.0);
\draw[black] (5.0, 2.0) -- (4.666667, 2.0);
\draw[black] (4.0, 7.5) to[out=0, in=180] (4.333333, 7.0);
\draw[black] (4.666667, 7.0) -- (4.333333, 7.0);
\draw[black] (5.0, 7.0) -- (4.666667, 7.0);
\draw[black] (4.0, 7.5) to[out=0, in=180] (4.333333, 8.0);
\draw[black] (4.666667, 8.0) -- (4.333333, 8.0);
\draw[black] (5.0, 8.0) -- (4.666667, 8.0);
\draw[black] (5.333333, 1.0) -- (5.0, 1.0);
\draw[black] (5.666667, 1.0) -- (5.333333, 1.0);
\draw[black] (6.0, 1.0) -- (5.666667, 1.0);
\draw[black] (5.333333, 2.0) -- (5.0, 2.0);
\draw[black] (5.666667, 2.0) -- (5.333333, 2.0);
\draw[black] (6.0, 2.0) -- (5.666667, 2.0);
\draw[black] (5.333333, 8.0) -- (5.0, 8.0);
\draw[black] (5.666667, 8.0) -- (5.333333, 8.0);
\draw[black] (6.0, 8.0) -- (5.666667, 8.0);
\draw[black] (5.333333, 7.0) -- (5.0, 7.0);
\draw[black] (5.666667, 7.0) -- (5.333333, 7.0);
\draw[black] (6.0, 7.0) -- (5.666667, 7.0);
\draw[black] (5.0, 4.5) to[out=270, in=180] (5.333333, 3.0);
\draw[black] (5.666667, 3.0) -- (5.333333, 3.0);
\draw[black] (6.0, 3.0) -- (5.666667, 3.0);
\draw[black] (5.0, 4.5) to[out=90, in=180] (5.333333, 6.0);
\draw[black] (5.666667, 6.0) -- (5.333333, 6.0);
\draw[black] (6.0, 4.0) -- (5.666667, 6.0);
\draw[black] (6.333333, 1.0) -- (6.0, 1.0);
\draw[black] (6.666667, 1.0) -- (6.333333, 1.0);
\draw[black] (7.0, 1.0) -- (6.666667, 1.0);
\draw[black] (6.333333, 2.0) -- (6.0, 2.0);
\draw[black] (6.666667, 2.0) -- (6.333333, 2.0);
\draw[black] (7.0, 2.0) -- (6.666667, 2.0);
\draw[black] (6.333333, 3.0) -- (6.0, 3.0);
\draw[black] (6.666667, 3.0) -- (6.333333, 3.0);
\draw[black] (7.0, 3.0) -- (6.666667, 3.0);
\draw[black] (6.333333, 4.0) -- (6.0, 4.0);
\draw[black] (7.0, 4.5) to[out=270, in=0] (6.666667, 4.0);
\draw[black] (6.666667, 4.0) -- (6.333333, 4.0);
\draw[black] (6.333333, 8.0) -- (6.0, 8.0);
\draw[black] (6.666667, 8.0) -- (6.333333, 8.0);
\draw[black] (7.0, 8.0) -- (6.666667, 8.0);
\draw[black] (6.333333, 7.0) -- (6.0, 7.0);
\draw[black] (6.666667, 7.0) -- (6.333333, 7.0);
\draw[black] (7.0, 7.0) -- (6.666667, 7.0);
\draw[black] (6.0, 5.5) to[out=270, in=180] (6.333333, 5.0);
\draw[black] (7.0, 4.5) to[out=90, in=0] (6.666667, 5.0);
\draw[black] (6.666667, 5.0) -- (6.333333, 5.0);
\draw[black] (6.0, 5.5) to[out=90, in=180] (6.333333, 6.0);
\draw[black] (6.666667, 6.0) -- (6.333333, 6.0);
\draw[black] (7.0, 6.0) -- (6.666667, 6.0);
\draw[black] (7.333333, 1.0) -- (7.0, 1.0);
\draw[black] (7.666667, 1.0) -- (7.333333, 1.0);
\draw[black] (8.0, 1.0) -- (7.666667, 1.0);
\draw[black] (7.333333, 2.0) -- (7.0, 2.0);
\draw[black] (7.666667, 2.0) -- (7.333333, 2.0);
\draw[black] (8.0, 2.0) -- (7.666667, 2.0);
\draw[black] (7.333333, 3.0) -- (7.0, 3.0);
\draw[black] (8.0, 4.5) to[out=270, in=0] (7.666667, 3.0);
\draw[black] (7.666667, 3.0) -- (7.333333, 3.0);
\draw[black] (7.333333, 8.0) -- (7.0, 8.0);
\draw[black] (7.666667, 8.0) -- (7.333333, 8.0);
\draw[black] (8.0, 8.0) -- (7.666667, 8.0);
\draw[black] (7.333333, 7.0) -- (7.0, 7.0);
\draw[black] (7.666667, 7.0) -- (7.333333, 7.0);
\draw[black] (8.0, 7.0) -- (7.666667, 7.0);
\draw[black] (7.333333, 6.0) -- (7.0, 6.0);
\draw[black] (8.0, 4.5) to[out=90, in=0] (7.666667, 6.0);
\draw[black] (7.666667, 6.0) -- (7.333333, 6.0);
\draw[black] (8.333333, 1.0) -- (8.0, 1.0);
\draw[black] (9.0, 1.5) to[out=270, in=0] (8.666667, 1.0);
\draw[black] (8.666667, 1.0) -- (8.333333, 1.0);
\draw[black] (8.333333, 2.0) -- (8.0, 2.0);
\draw[black] (9.0, 1.5) to[out=90, in=0] (8.666667, 2.0);
\draw[black] (8.666667, 2.0) -- (8.333333, 2.0);
\draw[black] (8.333333, 8.0) -- (8.0, 8.0);
\draw[black] (8.666667, 8.0) -- (8.333333, 8.0);
\draw[black] (9.0, 8.0) -- (8.666667, 8.0);
\draw[black] (8.333333, 7.0) -- (8.0, 7.0);
\draw[black] (8.666667, 7.0) -- (8.333333, 7.0);
\draw[black] (9.0, 7.0) -- (8.666667, 7.0);
\draw[black] (9.333333, 8.0) -- (9.0, 8.0);
\draw[black] (9.666667, 8.0) -- (9.333333, 8.0);
\draw[black] (10.0, 8.0) -- (9.666667, 8.0);
\draw[black] (9.333333, 7.0) -- (9.0, 7.0);
\draw[black] (9.666667, 7.0) -- (9.333333, 7.0);
\draw[black] (10.0, 7.0) -- (9.666667, 7.0);
\draw[black] (10.333333, 8.0) -- (10.0, 8.0);
\draw[black] (11.0, 7.5) to[out=90, in=0] (10.666667, 8.0);
\draw[black] (10.666667, 8.0) -- (10.333333, 8.0);
\draw[black] (10.333333, 7.0) -- (10.0, 7.0);
\draw[black] (11.0, 7.5) to[out=270, in=0] (10.666667, 7.0);
\draw[black] (10.666667, 7.0) -- (10.333333, 7.0);
\draw[black] (10.0, 1.5) to[out=270, in=180] (10.333333, 1.0);
\draw[black] (10.666667, 1.0) -- (10.333333, 1.0);
\draw[black] (11.0, 1.0) -- (10.666667, 1.0);
\draw[black] (10.0, 1.5) to[out=90, in=180] (10.333333, 2.0);
\draw[black] (10.666667, 2.0) -- (10.333333, 2.0);
\draw[black] (11.0, 2.0) -- (10.666667, 2.0);
\draw[black] (11.333333, 1.0) -- (11.0, 1.0);
\draw[black] (12.0, 1.5) to[out=180, in=0] (11.666667, 1.0);
\draw[black] (11.666667, 1.0) -- (11.333333, 1.0);
\draw[black] (11.333333, 2.0) -- (11.0, 2.0);
\draw[black] (12.0, 1.5) to[out=180, in=0] (11.666667, 2.0);
\draw[black] (11.666667, 2.0) -- (11.333333, 2.0);
\draw[black] (12.0, 1.5) to[out=0, in=180] (12.333333, 1.0);
\draw[black] (13.0, 1.5) to[out=270, in=0] (12.666667, 1.0);
\draw[black] (12.666667, 1.0) -- (12.333333, 1.0);
\draw[black] (12.0, 1.5) to[out=0, in=180] (12.333333, 2.0);
\draw[black] (13.0, 1.5) to[out=90, in=0] (12.666667, 2.0);
\draw[black] (12.666667, 2.0) -- (12.333333, 2.0);
\fill[red] (3.0, 1.5) circle (2.0pt);
\fill[red] (4.0, 7.5) circle (2.0pt);
\fill[red] (12.0, 1.5) circle (2.0pt);
\fill[blue] (1.0, 7.5) circle (2.0pt);
\fill[blue] (2.0, 1.5) circle (2.0pt);
\fill[blue] (5.0, 4.5) circle (2.0pt);
\fill[blue] (6.0, 5.5) circle (2.0pt);
\fill[blue] (7.0, 4.5) circle (2.0pt);
\fill[blue] (8.0, 4.5) circle (2.0pt);
\fill[blue] (9.0, 1.5) circle (2.0pt);
\fill[blue] (10.0, 1.5) circle (2.0pt);
\fill[blue] (11.0, 7.5) circle (2.0pt);
\fill[blue] (13.0, 1.5) circle (2.0pt);
   \end{tikzpicture}
   \caption{Topological structure of $Q(C_4)$: One oval and three figure eights with $Z_{11}$-singularities.
     In Fig.~\ref{fig:shustin_real} (left) the empty oval sits \enquote{vertically between} the two singularities on the left.}
   \label{fig:QC4}

   \vspace{1.4cm}
   \begin{tikzpicture}[scale=.245]
      \input{tikz/shustin_topological.tikz}
   \end{tikzpicture}
   \caption{Topological structure of $C_8$.
     One can clearly see four empty ovals and three ovals containing five empty ovals each.
     Thus the real scheme is $\scheme{4\sqcup 3\scheme{5}}$.}
   \label{fig:shustin_topological}
\end{figure}

Over time, many \mcurves of degree eight have been found, but historically the curve of Shustin marks a particularly challenging construction of a degree eight \mcurve with three non-empty ovals \cite{Shustin:1987}.
In the notation of Rokhlin--Viro the isotopy type of Shustin's curve reads $\scheme{4\sqcup 3\scheme{5}}$; that is, the curve consists of four empty ovals and three large ovals containing five empty ovals each.
Shustin's construction focuses on configurations and uses Viro's patchworking \cite{Viro:Patchworking}.
In this way, Shustin first constructs a degree eight singular curve, called $Q(C_4)$, which has three singularities of Arnold type $Z_{11}$ \cite{AGV}.
These singularities can be smoothed simultaneously to obtain the desired smooth curve, $C_8$.
We employ an observation of Gudkov \cite{gudkov1971} to modify Shustin's construction slightly, such that all intermediate steps remain rational.
In this way we obtain rational versions of $Q(C_4)$ and $C_8$; see Fig.~\ref{fig:shustin_real} for renderings of both curves.
As our first main contribution, we provide an explicit polynomial realizing Shustin's curve from \cite{Shustin:1987}, together with a computational proof in \oscar \cite{OSCAR,OSCAR-book} that the resulting curve satisfies all relevant properties.
We would like to point out that the isotopy type of Shustin's curve is inaccessible by combinatorial patchworking, that is, it is not realized by any T-curve.
This follows from a general argument \cite[Cor.~22]{121curves} based on Haas' zone decomposition \cite{1997:Haas:AlgebraicCombinatorial}.
However, those details are beyond our current scope.

The hardest part in our computational proof is the correct detection of the topological type of a real plane algebraic curve from its defining polynomial.
This can be achieved by quantifier elimination \cite{Collins:1975,Brown+Davenport:2007}, which works for real polynomials with any number of variables.
General implementations include \qepcad \cite{qepcadb}, \macaulay (open cells only) \cite{LeedelRioRahkooy:2503.21731}, and \mathematica.
In the plane, however, the algorithmic problem becomes simpler, allowing for code which is substantially faster.
In fact, detecting the topological type of a real affine plane curve is equivalent to drawing that curve in a topologically correct way; see Figs.~\ref{fig:QC4} and~\ref{fig:shustin_topological}.
Projective curves are handled by passing to suitable affine charts.
Algorithms for that task are standard, and previous implementations, such as \exacus \cite{exacus} or \cgal \cite{cgal}, mostly follow the same basic strategy:
First make the curve generic by means of a linear coordinate change so that no two points of interest are on the same vertical line.
Then identify all points with a vertical tangent, all singularities and maybe other points of interest.
Connect those points appropriately.
A notable deviation from this strategy is the algorithm of Diatta, Rouillier and Roy \cite{Diatta:2014}, which avoids any coordinate change; its refinement by Diatta, Diatta, Rouillier, Roy and Sagraloff \cite{Diatta:2022} matches the best known bit complexity $\tilde{O}(d^5\tau+d^6)$ for curves of degree~$d$ with integer coefficients of bitsize at most~$\tau$.
Working in the original coordinates preserves the sparsity of the input and yields a cylindrical algebraic decomposition of the curve itself, rather than of a sheared copy, at the price of handling vertical asymptotes and several critical points in the same vertical fiber.
A variant of this method, without a complexity analysis, is implemented in \curvetopology, which is part of the \julia package \pace \cite{pacejl}.
Nevertheless, many practical aspects and many design decisions remain nontrivial.
Our second main contribution is an implementation in \oscar, which is powerful enough to handle challenging input, including several smooth and singular curves arising in Shustin's construction \cite{Shustin:1987}.
We employ our implementation to verify the various steps of the construction.
For general input, the implementation checks its numerical decisions for consistency but does not certify them; see Section~\ref{sec:implementation}.
For the curves arising in Shustin's construction, however, all results have also been checked by hand, including the singular curve $Q(C_4)$ and its singularities, and the topology of $C_8$ was confirmed independently with \pace.

To the best of our knowledge, no explicit polynomial realizing Shustin's isotopy type has been published before; Shustin's construction proves the existence of suitable parameters but does not provide them.
Finding such a polynomial is not straightforward, since three $Z_{11}$-singularities have to be smoothed simultaneously and the resulting ovals are too small to be seen in a drawing with the actual coordinates.

This article is an expanded version of our extended abstract presented at MEGA 2026.
Since then, the drawing algorithm has been generalized to nongeneric input and to unbounded branches.
Previously, singular points were distinguished from vertical tangents by a heuristic; now the implementation relies on isolating boxes computed by \msolve.
Most importantly, the main example of the extended abstract was Gudkov's $M$-curve of degree six \cite{gudkov1971}; the construction of Shustin's curve of degree eight, including the verification of singularity types and final isotopy type, is entirely new.

The article is organized by reversing the order in the description above.
In Section~\ref{sec:drawing} we discuss drawing real plane curves.
Then, in Section~\ref{sec:implementation}, we explain the specifics of our \oscar implementation.
Finally, Section~\ref{sec:shustin} features a step-by-step construction of Shustin's curve, following the original article \cite{Shustin:1987} closely.
The code for our examples is available on GitHub at
\begin{center}
   \url{https://github.com/dmg-lab/DrawingCurvesInOscar}.
\end{center}

\section{Drawing curves} \label{sec:drawing}
Concerning the type of drawing of a plane curve, there is one fundamental decision to be made.
Either the goal is a topologically correct representation, i.e., the output is a planar graph which is isotopic to the real curve,
or the drawing tries to be visually correct, respecting ratios and the overall impression, but possibly lacking precision (due to lack of resolution) at special points.
Achieving both goals at the same time is often impossible, in particular, if the coordinates of the points of interest vary over several orders of magnitude.
Usually, one first computes a topologically correct representation; the state of the art is explained, e.g., in articles by Seidel and Wolpert \cite{Seidel+Wolpert:2005}, Cheng et al.\ \cite{ChengEtAl:2010,ChengEtAl:2024} and Diatta et al.\ \cite{Diatta:2014,Diatta:2022}.
Worst case complexity bounds have been given by Kerber and Sagraloff \cite{Kerber-Sagraloff:2012}.
Subsequently, one can try to draw the arcs, e.g., as B\'ezier curves.
Computing the points of interest of a curve, such as singularities, leads to solving systems of polynomial equations.
So the full array of techniques from resultants to Gröbner bases and numerical methods comes into play; cf.\ \cite{Sturmfels02-Solving}.

\subsection{Main algorithm}
We sketch our main procedure, which is a version of the method of Seidel and Wolpert \cite{Seidel+Wolpert:2005}.
Let $\field$ be some subfield of the reals.
For a bivariate polynomial $f\in \field[x,y]$, we want to develop a drawing of the real affine plane curve $V_\RR(f)$ in $\RR^2$, reflecting the topological properties of the curve. Note that the curve is not necessarily irreducible.
The input to our algorithm is always affine.
For projective curves, such as those in Section~\ref{sec:shustin}, we choose suitable affine charts by hand.

The algorithm requires generic input, i.e., the critical points of $f$ (see Step~(1) below) must have pairwise distinct $x$-coordinates, and the coefficient of $y^{\deg f}$ in $f$ must be nonzero.
This can always be achieved by a linear coordinate change, since a generic one works, and our implementation provides procedures for finding one; see below.
For generic $f$ the algorithm proceeds as follows.
\begin{enumerate}[(1)]
\item Solve the system $f=\partial f/\partial y=0$ over $\RR$ to get the \emph{critical points} of $f$ and sort them by $x$-coordinate.
  These are the points of $V_\RR(f)$ that are critical for the projection $(x,y)\mapsto x$, not the zeros of the gradient of $f$; each of them is either a singular point of $V_\RR(f)$ or a smooth point with a vertical tangent line.
  Since $f$ is generic, those $x$-coordinates are pairwise distinct.
\item For two subsequent critical points $(x_0, y_0)$ and $(x_3, y_3)$ set
  \[
    x_1:=x_0+1/3\cdot(x_3-x_0)\text{ and }x_2:=x_0+2/3\cdot(x_3-x_0).
  \]
\item Determine the real roots of the univariate polynomials $f_i:=f(x_i,y)$ for $0\leq i \leq 3$, and sort them in ascending order.
  We may take $x_1$ and $x_2$ rational; the critical fibers $f_0$ and $f_3$ are discussed in Section~\ref{sec:implementation}.
\item Roots $y_{i,j}$ of $f_i$ correspond to points $(x_i,y_{i,j})$ on $V_\RR(f)$.
  Now $f_1$ and $f_2$ have the same number of real roots.
  Indeed, the real roots of $f(t,y)$ depend continuously on $t\in(x_0,x_3)$, and their number could only change where two of them collide, i.e., at a critical point, or where one of them escapes to infinity, i.e., at a vertical asymptote.
  The former is excluded since there is no critical point with $x$-coordinate strictly between $x_0$ and $x_3$, and the latter since, by genericity, the coefficient of $y^{\deg f}$ in $f$ is a nonzero constant.
\item Exactly one root of $f_0$ and one root of $f_3$ correspond to critical points of $f$. Mark these.
\item Start by connecting the points corresponding to roots of $f_1$ and $f_2$ in their respective order.
  Since there are no critical points between $x_0$ and $x_3$, there cannot be any \enquote{crossings} or \enquote{turns}.
\item Now connect the points corresponding to unmarked roots of $f_0$ and roots of $f_1$.
  Connect the unmarked roots above the marked root of $f_0$ to the top roots of $f_1$.
  Proceed similarly from the bottom for those below.
\item Connect all remaining points corresponding to roots of $f_1$ to the critical point arising from the marked root of $f_0$.
  Do the same with $f_2$ and $f_3$.
\item To the left of the leftmost critical point $(x_0,y_0)$ pick $x_{-1}:=x_0-\delta$ for some $\delta>0$, determine the real roots of $f(x_{-1},y)$, and connect the corresponding points to those above $x_0$ as in Steps~(7) and~(8).
  Proceed similarly to the right of the rightmost critical point.
  As there are neither critical points nor vertical asymptotes beyond these two points, $V_\RR(f)$ is the disjoint union of graphs of finitely many continuous functions there, and the number of real roots of $f(t,y)$ does not depend on $t$.
  The resulting graph is isotopic to $V_\RR(f)$ with each unbounded branch truncated at $x_{-1}$ or at the corresponding point on the right.
\end{enumerate}

\begin{figure}[tb]
\includegraphics[width=.7\textwidth]{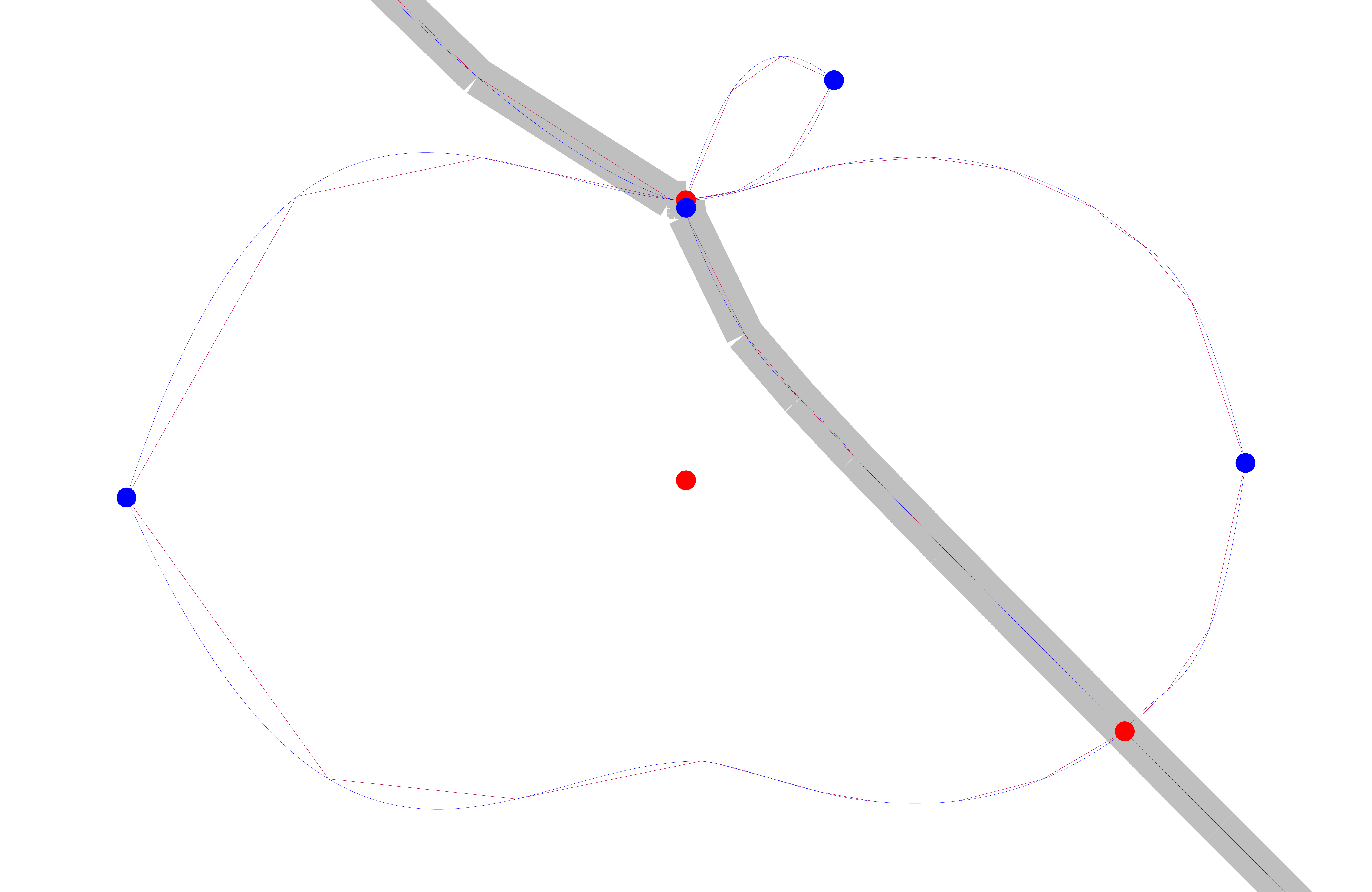}
\caption{The \enquote{apple} of \cite{Seidel+Wolpert:2005}.
  Piecewise-linear approximation in red and approximation by degree three B\'ezier curves in blue.
  \enquote{Diagonal} marked in gray for better comparison with Fig.~\ref{fig:apple_graph}.}
\label{fig:apple}
\end{figure}

Let $\overline{\field}$ denote the algebraic closure of $\field$.
The critical points and the roots of the $f_i$ lie in $(\overline{\field}\cap\RR)^2$ and $\overline{\field}\cap\RR$, respectively, where $\overline{\field}\cap\RR$ is the real closure of $\field$.
Hence the above procedure only involves arithmetic in $\field$ and in its real closure, and this arithmetic is crucial for the overall computational performance; we will give details in Section~\ref{sec:implementation} below.
The two nontrivial ingredients are closely related but distinct tasks:
\begin{enumerate}[(a)]
\item\label{it:critical} In Step~(1), find all critical points by solving a zero-dimensional system given by two bivariate polynomials.
\item\label{it:univariate} In Step~(3), find all real roots of the univariate polynomials $f_i$.
\end{enumerate}

In Step~(1) any number of further points of interest can be added to mark their correct positions; typical examples are inflection points.
We deliberately described our procedure with a rather vague notion of \enquote{drawing} as the output, as the specification of \enquote{drawing} depends on the purpose of the drawing.
As far as visual representations are concerned, there are at least two options: either exact piecewise-linear embeddings or piecewise-smooth embeddings via, e.g., B\'ezier curves; we will see both below.
Occasionally, it is useful to process the \enquote{drawing} further.
In that case, it may be useful to employ a suitable data structure to properly encode a (not necessarily connected) curve and the induced cell decomposition of $\RR^2$; e.g., doubly-connected edge lists \cite[\S2.2]{4marks}.
The resulting cell complex forms a \enquote{cylindrical algebraic decomposition} \cite{Collins:1975,Brown+Davenport:2007}; see the introduction for implementations in any number of variables.

Currently our implementation does not deal with polynomials without critical points.
This happens, e.g., if the real locus $V_\RR(f)$ is empty, or the curve consists of parallel lines.
We plan to revisit these situations in the future.
Our function \mintinline{jl}{draw_curve_tikz} returns a boolean.
A successful drawing is indicated by \mintinline{jl}{true}, otherwise \mintinline{jl}{false}.
To obtain generic input, \mintinline{jl}{draw_curve_tikz} applies linear coordinate changes until the transformed curve is generic, i.e., its critical points have pairwise distinct $x$-coordinates and the coefficient of $y^{\deg f}$ is nonzero.
Alternatively, the user may specify a transformation via the keyword argument \mintinline{jl}{transform}, together with a bound \mintinline{jl}{ntries} on the number of repetitions.
The isotopy graph is computed for the transformed curve, and its vertices are then mapped back to the original coordinates.
For Step~(9) we take $\delta$ to be one twentieth of the distance between the $x$-coordinates of the leftmost and the rightmost critical point, or $\delta=1$ if there is only one critical point.
An attempt with a given coordinate change also fails if the $x$-coordinates of the critical points cannot be separated at the chosen precision, which cannot be told apart from nongeneric input, or if $f_1$ and $f_2$ have different numbers of real roots, which by Step~(4) indicates a numerical failure.
If all attempts fail, \mintinline{jl}{false} is returned.
This may happen when a numerical backend for solving is used (see Section~\ref{sec:implementation} below); in that case one can try to increase the precision via the keyword argument \mintinline{jl}{solver_precision}.
For the computations in the present article we rely on exact arithmetic throughout.
We illustrate our implementation by taking as input the small example \enquote{apple} of \cite{Seidel+Wolpert:2005} (degree seven), which reads in \oscar code:

\begin{minted}{jl}
R, (x,y) = polynomial_ring(QQ, [:x,:y])
f = x^3*y^4 - 3*x^5 + 2*x^5*y^2 + 2*y^4 + x*y^4 - x^4 - 6*x^2*y^4 + x^3*y^2 + 4*x^2 + y^2 + 10*x^2*y^2 - 3*x^4*y^2 - 3*y^6 + x^7 - 4*x^3 + y^7 + 2*y^5 + y^3*x - y^2*x - y^5*x - 12*x^2*y + 4*x^3*y + 2*x^2*y^3 - 2*x^3*y^3 + 2*x^2*y^5 + 3*x^4*y - x^5*y + x^4*y^3 - 3*y^3
\end{minted}
The following features of the example stress-test any drawing algorithm:
\begin{enumerate}[(i)]
  \item an isolated singularity;
  \item a complicated singularity with six branches;
  \item several loops;
  \item a simpler singularity with four branches;
  \item critical points with the same $x$-coordinate.
\end{enumerate}
Moreover, the curve is not compact; its unbounded branches are handled by Step~(9).
Due to feature~(v), the apple is not generic in its original coordinates, and so the coordinate changes described above are essential.
The drawings in Figs.~\ref{fig:apple} and~\ref{fig:apple_graph} are produced by the following calls, writing \TikZ files; the option \mintinline{jl}{graph=true} selects the schematic rendering of Section~\ref{sec:topology}.
\begin{minted}{jlcon}
julia> draw_curve_tikz("apple_msolve.tikz", f)
true

julia> draw_curve_tikz("apple_msolve_graph.tikz", f; graph=true)
true
\end{minted}

\subsection{Schematic rendering}
\label{sec:topology}
Sometimes one is more interested in structural information of the curve, i.e., its topology, especially if the actual coordinates of the points from the algorithm lead to pictures that are hard to understand. Such a situation may arise if some points are very close together, while other points are several orders of magnitude further away. In this case, the points in close vicinity to each other will be hard to discern in the image or even print as one point due to lacking resolution.
If only the topology matters, however, one can choose sufficiently distinct coordinates for all points, but still connect them as in the original algorithm.

\begin{figure}
  \includegraphics[width=.55\textwidth]{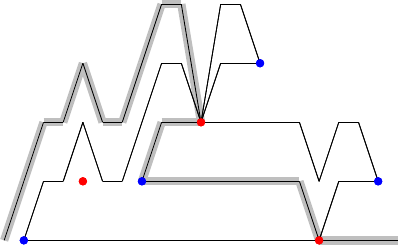}
  \caption{A schematic but topologically correct rendering of the \enquote{apple}; cf.\ Fig.~\ref{fig:apple}.
    Singular points are marked red, other critical points blue.
    The critical points refer to the transformed coordinates.}
\label{fig:apple_graph}
\end{figure}

In the algorithm, we first compute the critical points and sort them by $x$-coordinate.
In Step~(2) we choose some intermediate $x$-coordinates between two $x$-coordinates of adjacent critical points and denote the resulting sequence of coordinates on the $x$-axis 
by $\xi_1<\xi_2<\dots<\xi_n$. Subsequently, we compute the roots above every $\xi_i$, i.e., the real roots of $f(\xi_i, y)$; call them $\eta_{i,1}<\eta_{i,2}<\dots<\eta_{i,m_i}$.
(Note that the number $m_i$ obviously depends on $\xi_i$.)
Our schematic drawing uses the coordinates $(i,j)$ to encode the point $(\xi_i, \eta_{i,j})$.
Figure~\ref{fig:apple_graph} depicts our previous example schematically.

\section{The \oscar implementation}
\label{sec:implementation}
There are several viable ways to approach the crucial subproblems \eqref{it:critical} and \eqref{it:univariate} \cite{Sturmfels02-Solving}.
Two of them, elimination and homotopy continuation, are supported by \oscar \cite{OSCAR,OSCAR-book}.

\subsection{Elimination}

Computing a Gröbner basis of an ideal with respect to an elimination ordering is a standard method for solving a general polynomial system with finitely many solutions \cite{Sturmfels02-Solving}.
In the bivariate case, a reduced elimination Gröbner basis of a zero-dimensional ideal contains a univariate polynomial, the \emph{eliminant}.
Finding all (real) roots of the eliminant, substituting into the original system, and computing the (real) roots of the resulting univariate polynomials solves \eqref{it:critical}.
In this approach, the second subproblem \eqref{it:univariate} is necessary to compute all critical points.
The most general version of this approach is implemented in the \oscar command \mintinline{jl}{rational_solutions}.
However, it is often slow for relevant systems.

For us, a better alternative is \oscar's interface to \msolve \cite{berthomieu2021} via \algsolving.
First, it is usually faster.
Second, the command \mintinline{jl}{real_solutions} finds all solutions of the given system over $\RR$.
Assuming rational polynomials as input, the real solutions are, in fact, algebraic numbers; \algsolving employs interval arithmetic for representing such numbers.
Thus \msolve gives exactly one rational box per solution.

We continue to process the \enquote{apple}.
\begin{minted}{jlcon}
julia> critpt = ideal([f, derivative(f, y)]);

julia> rs,_ = real_solutions(Vector{Vector{QQFieldElem}}, critpt);

Restarting with another random linear form
Restarting with a non-random linear form

julia> rs
6-element Vector{Vector{Vector{QQFieldElem}}}:
 [[4880396824665781919//9223372036854775808, 9760793649331563913//18446744073709551616], [26193891127206274131//18446744073709551616, 26193891127206274329//18446744073709551616]]
 [[0, 0], [1, 1]]
 [[2, 2], [0, 0]]
 [[0, 0], [0, 0]]
 [[-2, -2], [0, 0]]
 [[14455604285078938641//9223372036854775808, 14455604285078939701//9223372036854775808], [-8269746749600746823//9223372036854775808, -16539493499201491451//18446744073709551616]]
\end{minted}
The output shows six two-dimensional boxes for the six critical points.
Some boxes, like \mintinline{jl}{[[0, 0], [1, 1]]}, shrink to single points.
In that example it means that the rational point $(0,1)$ is an exact solution.
For actually drawing the curves, the midpoints of the boxes serve as (overly precise) approximations to the nonrational solutions.

For the critical fibers $f_0$ and $f_3$ of Step~(3), \msolve only provides a box $[a,b]\times[c,d]$ around the critical point.
If $a<b$, we replace the critical fiber by the regular fiber over $a$ or $b$ with more real roots, and we merge its roots in $[c-\epsilon,d+\epsilon]$ into the critical point.
The tolerance $\epsilon$ is chosen adaptively.
Starting from $\epsilon=2^{-k}$, for $k$ given by the keyword argument \mintinline{jl}{ybox_tolerance}, we halve $\epsilon$ until three conditions hold:
no root of the fibers over $a$ and $b$ lies in $[c-2\epsilon,d+2\epsilon]$ but outside $[c-\epsilon,d+\epsilon]$;
the two fibers have the same number of roots below this window, and the same number above it;
and, at a vertical tangent, the window contains two roots of one fiber and none of the other.
If no $\epsilon\geq 2^{-\lfloor p/2\rfloor}$ satisfies these conditions, where $p$ is given by \mintinline{jl}{solver_precision}, the attempt fails.
This is a consistency check rather than a proof, but it replaces a fixed guess for $\epsilon$.

To decide whether a critical point is singular or smooth with a vertical tangent, we also compute the singular points, i.e., the solutions of $f=\partial f/\partial x=\partial f/\partial y=0$, with \msolve.
A critical point is singular if and only if its box meets the box of a singular point; we check that each critical box meets at most one singular box, and that each singular point is found in this way.
This test is rigorous, since \msolve returns isolating boxes.

\subsection{Homotopy continuation}

Homotopy continuation is a method which fuses algebraic techniques with numerical analysis.
This approach is based on the following idea.
One starts with a polynomial system whose solutions are known.
That start system is continuously deformed into the given system, and the known solutions are numerically traced to the desired ones.

A sophisticated implementation is \homcon \cite{breiding2018}.\footnote{\homcon can be used as a package extension to \oscar \cite{OSCAR,OSCAR-book}.}
However, solutions obtained in this way are usually not suitable for further symbolic computations.
More crucially, the precision is limited by \julia's \mintinline{jl}{Float64}, which is used by \homcon.
This is not always enough to distinguish critical points from each other in case they are close together, while other critical points are far apart.
Still, \homcon is exceptionally fast compared to symbolic methods and delivers correct pictures for many example systems.

\subsection{Coefficient fields}
There are various solvers in and around \oscar one can use as a backend.
In all cases, drawing the curve with its actual coordinates eventually requires a conversion to floating point numbers.
This loss of precision may lead to collisions of distinct points, and then the picture no longer shows the topology correctly.
The schematic rendering of Section~\ref{sec:topology} avoids this problem: it only uses the order of the coordinates, which, with \msolve as the backend, is read off from exact isolating intervals before any conversion takes place.

For exact computations, \oscar offers a variety of coefficient fields that the user may choose from.
The most basic field is \mintinline{jl}{QQ}, the field of rational numbers, based on the \flint~\cite{flint} implementation. This is the only field that \msolve currently supports for our drawing algorithm.

Apart from \mintinline{jl}{QQ}, there are other relevant exact fields, like \mintinline{jl}{algebraic_closure(QQ)}, realized by \calcium, which is now merged into \flint.
Then, via \hecke \cite{nemo}, one can use \mintinline{jl}{number_field}.
In principle, \oscar provides solutions for polynomial systems over these fields, via the function \mintinline{jl}{rational_solutions}.
However, as these fields get increasingly complicated, performance becomes a limiting factor.
For more complicated input, as in Section~\ref{sec:shustin} below, this approach is not feasible.
Besides that, \oscar also has the inexact fields \mintinline{jl}{RR} and \mintinline{jl}{CC}. Both can be used in conjunction with \homcon, if one accepts the conversion to \mintinline{jl}{Float64}.

\section{Shustin's \mcurve in \oscar}
\label{sec:shustin}

This section now covers the example whose study caused us to embark on a journey to draw topologically correct images of plane curves. 
We use \oscar to follow Shustin's multi-step procedure \cite{Shustin:1987} in an explicit computation, obtaining polynomials with rational coefficients for the relevant curves. Our schematic drawings allow us to verify the topology, which was presented in the original article by means of elaborate hand-drawn sketches.

Shustin's construction is highly abstract and follows the classical philosophy which grew out of Gudkov's work: First construct a curve of the desired degree with carefully chosen and carefully placed singularities, then use a carefully chosen smoothing to obtain the desired topological structure. The subtleties of the construction are in the choices. For the choice of the singularity, the so-called classifier of Arnold \cite{AGV} is of crucial importance, as it provides information both on the initial form of the power series describing a given singularity and on the topology of a smoothing thereof. For the choice of the smoothing, it is vital to observe that the deformation may be performed using a curve which again contributes its own topology. It is the particular strength of Viro's method that it allows singularities of the chosen curve to be smoothed independently of each other.
This is vital to have sufficient control over the relative positioning of the ovals.

\begin{figure}[th]
   \begin{tabular}{cc}
      \includegraphics[height=4.5cm]{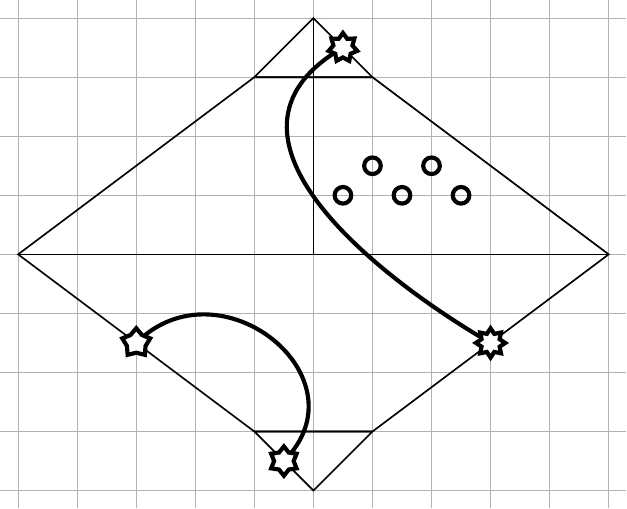}
      &
      \includegraphics[height=4.5cm]{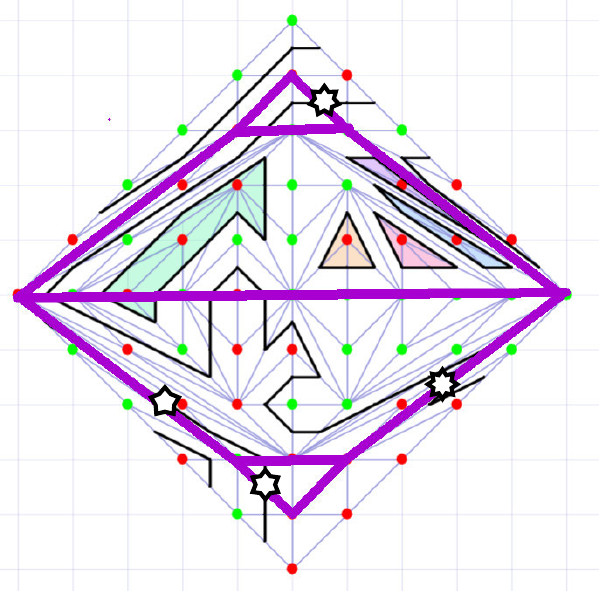}
   \end{tabular}
   \caption{Shustin's Fig.~1 and our patchworked $C_5$.
     Most relevant for Viro's patchwork is that the curves agree on their charts.
     In our case this means that the curve intersects the correct line segments, here marked as stars.
     On the right, the boundaries of the Newton polytopes of $C_5$ and $\ell$, reflected to all four quadrants, are drawn in purple.}
   \label{fig:ShustinMarcel}
\end{figure}

We will use Shustin's notation to be consistent with the original article.
First one constructs a curve $Q(C_4)$ with three $Z_{11}$-singularities in Arnold's notation.
Graphically speaking, $Q(C_4)$ consists of one empty oval and three figure eight components, with said singularities at their crossing points, see Fig.~\ref{fig:QC4}.
These singularities are then smoothed simultaneously using Viro's patchworking method relying in particular on a smooth quintic with prescribed topology from work of Gudkov \cite[Fig.~11]{Gudkov:1974}.
By doing this, one half of the figure eight becomes an empty oval, and the other half becomes an oval with five small ovals inside; the part that gets inserted is depicted in Fig.~\ref{fig:ShustinMarcel}.
At this point, the entire procedure can easily be sketched by hand.
However, it is much harder to find an actual real bivariate polynomial realizing this curve.
Indeed, Shustin never provides such a polynomial; instead he proves that suitable parameters exist and that smoothing the singularities simultaneously is possible.

Thus, we set out to construct the homogeneous polynomial of degree eight describing this curve in \oscar.
All curves in the construction are projective; whenever we analyze singularities or draw a curve, we pass to an affine chart suited to the step at hand.
There are two main challenges:
First, at any step, we want to verify correctness, that is the topological structure and the singularity types. Obtaining a topologically correct visual representation of the curve will be essential at this step.
Second, performance degrades dramatically if the coefficients are irrational real numbers.
Due to smoothness, it is already clear that the final polynomial may be chosen with coefficients in $\QQ$. We can just add a tiny perturbation to every coefficient to make it rational, and we will keep both smoothness and topological structure.
However, an elegant trick hidden in a paper of Gudkov \cite{gudkov1971} allows us to slightly modify Shustin's construction and keep all intermediate steps rational.

The final step of smoothing all three $Z_{11}$-singularities is depicted in Fig.~\ref{fig:shustin_real}.
The changes to the topological structure are barely visible to the naked eye; in particular, the inserted ovals cannot be detected.
Thus, Fig.~\ref{fig:shustin_topological} depicts the topological structure of the final curve, and, after the upcoming discussion, even serves as a proof for the claimed isotopy type.

\subsection{Overview of the construction}
Let us start by summarizing Shustin's construction from \cite{Shustin:1987}.

\begin{enumerate}[(1)]
\item Start from a suitable quartic:
Take two ellipses $C_1$ and $C_2$ that intersect in four points.\footnote{We follow Gudkov's and Shustin's abuse of notation denoting both the curve and its defining polynomial by the same label.}
Now perturb the defining polynomial $C_1 \cdot C_2$ of their union to obtain a curve $C_4$ consisting of four ovals.
\item Choose inflection points and tangents:
Each of the four ovals arises from one of the four \emph{lunes} between $C_1$ and $C_2$, i.e., the regions inside exactly one of the two ellipses.
The boundary of such a lune consists of an arc of $C_1$ and an arc of $C_2$, which bend in opposite directions as seen from the lune.
Hence, for a small perturbation, each oval has two inflection points, one near each corner of its lune.
Pick three of these inflection points $p_0$, $p_1$, and $p_2$ on three different ovals and let $L_i$ be the tangent line to the respective oval at $p_i$.
Choose suitable coordinates such that the $L_i$ become the coordinate axes.
(This requires the three tangents to be non-concurrent; see Section~\ref{sec:inflection} below for our explicit choices.)
\item Apply a Cremona transformation:
Now apply the transformation $Q: \PP^2 \dashrightarrow \PP^2$ taking $(x:y:z)$ to $(yz : xz : xy)$; see any textbook on algebraic geometry, e.g., \cite{Sha,Har}.
One obtains a degree eight curve $Q(C_4)$ consisting of one oval and three figure eights, where the three singular points are stated to be $Z_{11}$ singularities.
\item Smooth the curve:
By Viro's patchworking construction the defining polynomial is now changed in a very controlled way, stated as adding a quintic and a line in each of the three affine charts. Through this construction the curve then becomes smooth and has the expected topology.
\end{enumerate}

We now retrace this construction step by step in the following subsections.

\subsection{Prescribing inflection points}
\label{sec:inflection}
It is easy to choose $C_1$ and $C_2$ such that they intersect in four rational points.
As explained in step~(2) of the overview, a small perturbation yields four ovals with two inflection points each.
However, these inflection points tend not to be rational.

While working through the references of Shustin's work \cite{Shustin:1987}, we also went
through work of Gudkov \cite{gudkov1971}. The procedure there is very similar;
however, Gudkov perturbs the product $C_1\cdot C_2$ in a more subtle way. He
fixes a line $\ell$ that does not meet $C_1$ or $C_2$ in any real point, and the tangent $L_1$ to
$C_1$ at an intersection point $p_1$ of $C_1$ and $C_2$. Then he takes the
perturbation
\[
   C_4\ :=\ C_1\cdot C_2 + \lambda\cdot L_1\cdot \ell^3.
\]
Now Gudkov states that $L_1$ is still tangent to $C_4$ in $p_1$, and that $p_1$ has become an inflection point of $C_4$.

Our idea is to extend that trick and apply it thrice in our situation to
three of the intersection points of $C_1$ and $C_2$.
Since the intersection points were already rational, this will keep the
inflection points rational, too.

\begin{proposition}\label{prop:inflection}
  Let $C_1$ and $C_2$ be as above and take their union.
  Now take three of the intersection points $p_0$, $p_1$, and $p_2$.
  At every $p_i$ let $L_i$ be the tangent to one of the two components.
  Then take the perturbation
  \[
    C_4\ :=\ C_1\cdot C_2 + \lambda\cdot L_0\cdot L_1\cdot L_2\cdot \ell,
  \]
  where $\lambda\neq0$, and $\ell$ is a line which does not meet $C_1$ or $C_2$ in any real point.

  Then the points $p_i$ still lie on $C_4$ and the $L_i$ are still tangent to $C_4$ at the $p_i$.
  Furthermore the $p_i$ have become inflection points of $C_4$.
\end{proposition}

\begin{proof}
Without loss of generality we only prove the statement for $L_0$ tangent to $C_1$ at $p_0$, as the other statements follow by symmetry.

First observe that $p_0$ still lies on $C_4$, that is the defining polynomial $C_4$ vanishes at $p_0$, since
\[
   C_1(p_0)\ =\ C_2(p_0)\ =\ L_0(p_0)\ =\ 0.
\]
From the differentials, it is nearly immediate that $L_0$ is tangent to $C_4$ at $p_0$:
\[
   dC_4\ =\ C_1\cdot dC_2+C_2\cdot dC_1+\lambda\cdot (L_0\cdot d(L_1\cdot L_2\cdot \ell)+L_1\cdot L_2\cdot \ell\cdot dL_0).
\]
Evaluating at $p_0$ this becomes
\[
   dC_4(p_0) = \lambda\cdot L_1(p_0)\cdot L_2(p_0)\cdot dL_0(p_0).
\]
So $dC_4(p_0)$ is just a constant multiple of $dL_0(p_0)$. Moreover, $L_1(p_0)\cdot L_2(p_0) \neq 0$, because by B\'ezout and the tangency of $L_i$ to one ellipse at another intersection point $p_i$ no further intersection point $p_j$ lies on the line $L_i$ ($j \neq i$). Hence the tangent to $C_4$ at $p_0$ is $L_0$. 

There are several possibilities for showing that $p_0$ is an inflection point of $C_4$.
Computing via the Hessian is possible, but tedious, so instead we will show that the intersection multiplicity of $C_4$ and $L_0$ at $p_0$ is three.
We get
\[
   \begin{array}{ccl}
   I_{p_0}(C_4,\ L_0) &=& I_{p_0}(C_1\cdot C_2\ +\ \lambda\cdot L_0\cdot L_1\cdot L_2\cdot \ell,\ L_0)\\
   &=& I_{p_0}(C_1\cdot C_2,\ L_0)\\
   &=& I_{p_0}(C_1,\ L_0)+I_{p_0}(C_2,\ L_0).
   \end{array}
\]

By assumption, $C_1$ and $C_2$ meet in four distinct points, hence, by B\'ezout, they meet transversally in every point.
Thus, since $L_0$ is tangent to $C_1$ at $p_0$, it is not tangent to $C_2$ at $p_0$ and instead intersects $C_2$ transversally.
Tangency of $L_0$ to $C_1$ at $p_0$ gives $I_{p_0}(C_1, L_0)=2$, transversality gives $I_{p_0}(C_2, L_0)=1$, whence we indeed obtain $I_{p_0}(C_4, L_0)=3$.
\end{proof}

Proposition~\ref{prop:inflection} is local, and it does not say anything about the global shape of $C_4$.
For our explicit construction we take the circle $C_1=\tfrac{1}{25}(x^2+y^2)-z^2$ and the ellipse $C_2=(\tfrac{15}{53})^2x^2+(\tfrac{7}{53})^2y^2-z^2$, which meet in the four rational points $(\pm3:\pm4:1)$.
We pick $p_0=(3:4:1)$ with $L_0$ tangent to $C_2$, and $p_1=(-3:4:1)$, $p_2=(-3:-4:1)$ with $L_1$, $L_2$ tangent to $C_1$; moreover, $\ell=x-8z$ and $\lambda=\tfrac{1}{512}$.
That $C_4$ consists of four ovals with $p_0$, $p_1$, $p_2$ on three different ones, that $L_0$, $L_1$, $L_2$ are not concurrent, and that $C_4$ avoids the base points $L_i\cap L_j$ of the Cremona transformation is verified by a script in our code repository.
Applying the Cremona transformation to the resulting curve $C_4$ yields Shustin's $Q(C_4)$; see also Fig.~\ref{fig:QC4}.

\subsection{Verifying the singularity types}
The goal of this section is to verify that our $Q(C_4)$ has the singularities claimed in Shustin's article \cite{Shustin:1987}.
\begin{proposition}\label{prop:Z11}
The degree eight curve $Q(C_4)$ has exactly three singular points, and each of them is of type $Z_{11}$ in Arnold's classification.
\end{proposition}
   
Using \oscar we first verify that the singular locus really consists of three points, by computing the locus of simultaneous vanishing of the equation and its partial derivatives: 

\begin{minted}{jlcon}
julia> QC4 = load("QC4.mrdi");

julia> Rh = parent(QC4)
Multivariate polynomial ring in 3 variables over QQ graded by
  x -> [1]
  y -> [1]
  z -> [1]

julia> sing_ideal = ideal(QC4) + jacobi_ideal(QC4);

julia> rational_solutions(sing_ideal)
3-element Vector{Vector{QQFieldElem}}:
 [0, 1, 0]
 [1, 0, 0]
 [0, 0, 1]
\end{minted}

Arnold's classifier \cite{AGV} determines the type of a singularity and a normal form thereof up to analytic coordinate change (for over 100 cases with sufficiently low invariants). By this classifier, a $Z_{11}$ singularity can be recognized by the following data of a power series $f$ locally describing the singularity in question: $\operatorname{ord}(f)=4$, the initial part consists of a triple and a single factor and the Milnor number is $11$. The order of the power series and the presence of a triple and a single factor in the degree four part are rather obvious sorting criteria. They can be verified directly using \oscar:

\begin{minted}{jlcon}
julia> QC4 = load("QC4.mrdi");

julia> Rh=parent(QC4);                  # base_ring for projective curve

julia> R,(x,y) = polynomial_ring(QQ,[:x, :y]);  

julia> zchart = hom(Rh, R, [x,y,1]);      # map to affine Chart

julia> zQC4 = zchart(QC4);                # curve in affine Chart

julia> normal_form(zQC4,ideal(R,[x,y])^4) # no terms of degree less than 4
0

julia> j4QC4 = normal_form(zQC4,ideal(R,[x,y])^5)  # all terms up to degree 4

2359824029541015625//186524912152562171904*x^4 - 46714976316671875//971483917461261312*x^3*y + 7435875868735625//121435489682657664*x^2*y^2 - 370638183178825//11384577157749156*x*y^3 + 1972628140481//316238254381921*y^4

julia> factor(j4QC4)                      # factorize
(7890481//186524912152562171904) * (25*x - 16*y)^3 * (19140625*x - 36000144*y)
\end{minted}

The Milnor number is a more subtle measure for the complexity of the singularity. Roughly speaking, the Milnor number $\mu$ counts the spheres arising in a smoothing of the singularity by adding a nonzero constant. For an isolated hypersurface singularity at the origin of the affine plane, it can be computed as 
\[\mu (f) = \dim_\CC \CC\{x,y\}/ \big\langle \frac{\partial f}{\partial x},
                                      \frac{\partial f}{\partial y} \big\rangle.\]
Hence the verification boils down to a simple standard basis computation with respect to a local monomial ordering (e.g., negdeglex)---also available in \oscar via the command \mintinline{jl}{milnor_number}. In the subsequent code snippet the passage to the local analytic setting is performed explicitly, but without requiring knowledge of standard bases:

\begin{minted}{jlcon}
julia> QC4 = load("QC4.mrdi");

julia> Rh=parent(QC4);

julia> R,(x,y) = polynomial_ring(QQ,[:x, :y]);

julia> zchart = hom(Rh, R, [x,y,1]);

julia> zQC4 = zchart(QC4);

julia> R_curve =quo(R,zQC4)[1]                  # ring of the affine curve
Quotient
  of multivariate polynomial ring in 2 variables x, y
    over rational field
  by ideal with 1 generator

julia> # pass to analytic setting by forming the germ at the origin

julia> zQC4_germ = hypersurface_germ(R_curve,[0,0])[1]

julia> milnor_number(zQC4_germ)
11
\end{minted}

We only presented the $z$-chart in the above snippets; the computations in the other two charts are analogous, and the script in our code repository computes the Milnor numbers in all three charts.
This proves Proposition~\ref{prop:Z11}.

\subsection{A suitable transformation}
Classifying the singularity at the origin of a chart as $Z_{11}$ according to Arnold's list implies that there is an analytic coordinate change which transforms it into the so-called normal form
\[ x^5 - xy^3.\]
For our purposes, it suffices to find a suitable affine transformation which will turn the Newton boundary of $Q(C_4)$ into the edge segments $[(0,4),(1,3)]$ and $[(1,3),(5,0)]$.
In practice, one uses the factorization of the 4-jet of $Q(C_4)$, as computed above, to find its triple factor, which was $(25\cdot x -  16\cdot y)^3$.
A change of coordinates such that the $x$-terms cancel then yields the desired Newton boundary.

\begin{minted}{jlcon}
julia> vertices(newton_polytope(zQC4))
3-element SubObjectIterator{PointVector{QQFieldElem}}:
 [4, 4]
 [4, 0]
 [0, 4]

julia> zQC4mod = zQC4(16*x+100*y, 25*x+y);

julia> vertices(newton_polytope(zQC4mod))
5-element SubObjectIterator{PointVector{QQFieldElem}}:
 [8, 0]
 [5, 0]
 [1, 3]
 [0, 8]
 [0, 4]
\end{minted}

For the sake of a better understanding of the next steps, we include a monomial diagram of the situation after coordinate change, see Fig.~\ref{fig:monomial}. Changes to monomials in the darker shaded area will under no circumstances change the type of the singularity according to a determinacy result also used in the classifier. Moreover, the only monomial of degree five relevant for the type is the pure power in $x$. Passing to our curve in a different affine chart, this implies that we have quite some freedom to add low degree monomials there without affecting the singularity type in the previous chart. This is vital for understanding the independence in the next subsection.

\begin{figure}
   \begin{tikzpicture}[scale=.5]
      \fill[black!40] (8,0) -- (5.5,0) -- (0,5.5) -- (0,8) -- (8,8) -- cycle;
      \fill[black!20] (5.5,0) -- (5,0) -- (1,3) -- (0,4) -- (0,5.5) -- cycle;
      \draw[thin, black!30] (-.3,-.3) grid (8.3,8.3);
      \fill[black] (5,0) node[anchor=north east] {$x^5$} circle (4pt);
      \fill[black] (1,3) node[anchor=north east] {$xy^3$} circle (4pt);
      \fill[black] (0,4) node[anchor=north east] {$y^4$} circle (4pt);
      \draw[thick] (5,0) -- (1,3);
      \draw[thick] (0,4) -- (1,3);
      \draw[thick, dashed] (0,5.5) -- (5.5,0);
   \end{tikzpicture}
   \caption{Monomial diagram of a $Z_{11}$ singularity. Changes to the darker shaded area cannot change the singularity type because a singularity of this type is $5$-determined. According to more subtle arguments in the proof of Arnold's classifier, also the monomials in the light gray area do not change the type. This gives an impression of the freedom we have for changing lower degree monomials in the other charts without destroying this $Z_{11}$.}
   \label{fig:monomial}
\end{figure}
Please note that the coordinate change is applied locally in a single chart, which is only possible due to the independence of the Viro patchwork steps in the individual charts. For recombination at the end, the coordinate changes need to be reverted before gluing the charts back together. 

\subsection{Viro's patchwork}
We do not explain Viro's patchwork in full generality, for details please see \cite{Viro:Patchworking}.
Instead we will be very specific to our setting.
Having brought the respective chart of $Q(C_4)$ into the desired shape,
Shustin now wants to obtain a very particular smoothing by adding two further polynomials, a quintic $C_5$ and a line $\ell$.
Together the Newton polytopes of $Q(C_4)$, $C_5$, and $\ell$ introduce a regular subdivision of $8\cdot\Delta^2$, as depicted in Fig.~\ref{fig:NewtonSubdivision}.
A height function $8\cdot\Delta^2\to\RR$ inducing this regular subdivision is
\[
   \max\{0, -3\cdot x-4\cdot y+15, -3\cdot x-3\cdot y+12\}.
\]

\begin{figure}
   \begin{tikzpicture}[scale=.5]
      \fill[black!20] (8,0) -- (5,0) -- (1,3) -- (0,4) -- (0,8) -- cycle;
      \fill[black!30] (0,3) -- (1,3) -- (0,4) -- cycle;
      \fill[black!40] (0,3) -- (1,3) -- (5,0) -- (0,0) -- cycle;
      \draw[thin, black!50] (-.3,-.3) grid (8.3,8.3);
      \node[fill=white, anchor=south west] at (4,4) {$Q(C_4)$};
      \node[fill=white, anchor=east] at (0,1) {$C_5$};
      \node[fill=white, anchor=east] at (0,3.3) {$\ell$};
   \end{tikzpicture}
   \caption{Subdivision induced by the Newton polytopes of $Q(C_4)$, $C_5$, and $\ell$.}
   \label{fig:NewtonSubdivision}
\end{figure}

Shustin refers to \cite{Gudkov:1974} for the construction of $C_5$.
However, in degree five every isotopy type can be achieved by combinatorial patchwork.
Thus we can use a tool like
\begin{center}
   \url{https://patchworking.polymake.org/}
\end{center}
to achieve the desired type, see also \cite{121curves}.

In Fig.~\ref{fig:ShustinMarcel} on the left we have reproduced Fig.~1 of \cite{Shustin:1987}.
Note that the exact polynomials describing the curves $C_5$ and $\ell$ do not matter, as long as the boundaries of their Newton polytopes are compatible with the setting in the following sense:
\begin{itemize}\raggedright
   \item No intersection with the line segments $[(1,3),(5,0)]$, $[(-1,3),(-5,0)]$, $[(-1,3),(0,4)]$, and $[(0,-4),(1,-3)]$.
   \item Exactly one intersection with the line segments $[(1,3),(0,4)]$, $[(1,-3),(5,0)]$, $[(0,-4),(-1,-3)]$, and $[(-1,-3),(-5,0)]$.
\end{itemize}
On the right hand side of Fig.~\ref{fig:ShustinMarcel} we have depicted the combinatorial patchwork of our $C_5$ and $\ell$ produced with the aforementioned tool.
Looking at the boundary added in purple, we can verify that it satisfies exactly the behavior we desired.

The last constraint on the choice of $C_5$ and $\ell$ is to make sure that the polynomials $Q(C_4)$, $C_5$, and $\ell$ agree on common boundaries of their Newton polytopes.
That means that the coefficients in front of \enquote{common} monomials must be the same, e.g., $Q(C_4)$ and $C_5$ must have the same coefficients in front of $x\cdot y^3$ and $x^5$.
There are four such monomials to be considered: $y^3$, $y^4$, $x\cdot y^3$, and $x^5$.
Each of the three pairs of polynomials shares exactly two of them, since the common edges $[(0,3),(1,3)]$, $[(1,3),(0,4)]$, and $[(1,3),(5,0)]$ contain no further lattice points.
Multiplication with a nonzero scalar and rescaling of variables does not change the isotopy type of a polynomial.
So we introduce, for each of the three polynomials, a scalar factor and scaling factors for $x$ and $y$; equating the coefficients of the shared monomials yields six polynomial equations in these nine unknowns, which we solve using \oscar.
We are left with modified $Q(C_4)$, $C_5$, and $\ell$ that satisfy the initial conditions for Viro's patchwork.

In practice, we have written three methods in \oscar:
\begin{itemize}
   \item \mintinline{jl}{_patchworks2polynomial}: This function takes a file saved by the combinatorial patchworking tool and turns it into an \oscar polynomial.
   \item \mintinline{jl}{adjust_poly_overlaps}: This function takes a list of polynomials and solves a polynomial system to make the polynomials agree on the overlaps of their Newton polytopes.
      It returns a list of adjusted polynomials that still have the same curve isotopy types.
   \item \mintinline{jl}{construct_patchworked_poly}: This function takes a list of polynomials that satisfy Viro's patchworking conditions and a convex function inducing a regular subdivision into the Newton polytopes of these polynomials and constructs the resulting patchworked polynomial.
\end{itemize}
Warning: As of this writing, none of the above three functions does any input sanity checks. They were programmed for this particular example and would need significant refinement for any other application. They may, however, serve as code examples for similar tasks, whence we mention them.

In the following code snippets the polynomials \mintinline{jl}{ztriangleP}, \mintinline{jl}{zquadrangleP}, and \mintinline{jl}{zQC4mod} are patched together.
The first two polynomials belong to the triangle $\conv\{[0,3], [1,3], [0,4]\}$ and quadrangle $\conv\{[0,0],[5,0],[1,3],[0,3]\}$ respectively, and have the properties required by Shustin.

\begin{minted}{jlcon}
julia> ztriangleP = zchart(triangleP)
x*y^3 - 1//58149737003040059690390169*y^4 + 2954312706550833698643*y^3

julia> zquadrangleP = zchart(quadrangleP);

julia> vertices(newton_polytope(ztriangleP))
3-element SubObjectIterator{PointVector{QQFieldElem}}:
 [1, 3]
 [0, 4]
 [0, 3]

julia> vertices(newton_polytope(zquadrangleP))
4-element SubObjectIterator{PointVector{QQFieldElem}}:
 [5, 0]
 [1, 3]
 [0, 3]
 [0, 0]

julia> coeff(zQC4mod, [1,3]) == coeff(zquadrangleP, [1,3])
false

julia> polys, transform = adjust_poly_overlaps([zQC4mod, ztriangleP, zquadrangleP]);

julia> coeff(polys[1], [5,0])
-172539992867247277185821771427030500642102380884587871226791338776115529154632329345773282761727812180611405951060510270813456644283072880400998825968906034842465292503171540969209429928590149250138280282511690070537833570272195757786191555725132368429671048636182444773584142396962362806873438196908724538126202569680562134422629142501//13818940737753730222585759523950022304238617971870243881598491340410700127185738411400910700615554826240000000000

julia> coeff(polys[1], [5,0]) == coeff(polys[3], [5,0])
true

julia> coeff(polys[1], [1,3]) == coeff(polys[3], [1,3])
true

julia> coeff(polys[2], [1,3]) == coeff(polys[3], [1,3])
true

julia> tilt = matrix(ZZ, [0 0 0; -3 -4 15; -3 -3 12])
[ 0    0    0]
[-3   -4   15]
[-3   -3   12]

julia> pwp = construct_patchworked_poly(Rt, tilt, polys)
\end{minted}

\subsection{Loading from file and verification}
We have stored the final curve in the file \mintinline{bash}{shustin_C8.mrdi} that can be downloaded from \cite{zenodo:shustin}.
It is stored in the \mintinline{bash}{mrdi} format \cite{DellaVecchiaJoswigLorenz:2024}.
It can be loaded in \oscar and one can use the drawing algorithm described above for visualizing the topological structure.
As this curve is smooth, all its critical points are vertical tangents; there, the consistency check for $\epsilon$ in Section~\ref{sec:implementation} is particularly restrictive, since exactly two roots must enter the window, from one side only.
While $C_8$ lives in $\RR\PP^2$, our drawing algorithm takes affine input; so we need to pick affine charts.
In general, several charts may be necessary to get a complete view of the topological structure.
For $C_8$, however, a single chart suffices.
After the linear coordinate change $B$ in the session below, the binary form $C_8(x,y,0)$ has no real zeros, i.e., the line at infinity $\{z=0\}$ does not meet the real curve.
Hence $V_\RR(C_8)$ is entirely contained in the affine chart $\{z\neq0\}\cong\RR^2$, and the mutual position of the ovals in that chart determines the isotopy type in $\RR\PP^2$.
Such a chart can only exist for curves of even degree, as a real curve of odd degree meets every real line.

\begin{theorem}
The curve $C_8$ from our file \mintinline{bash}{shustin_C8.mrdi} has the isotopy type claimed by Shustin; it is an \mcurve of type $\scheme{4\sqcup 3\scheme{5}}$.
\end{theorem}

\begin{minted}{jlcon}
julia> C8 = load("shustin_C8.mrdi");

julia> B = transpose(matrix(QQ, ([1 0 1; 0 1 1; 0 0 1])))
[1   0   0]
[0   1   0]
[1   1   1]

julia> Rh = parent(C8)
Multivariate polynomial ring in 3 variables over QQ graded by
  x -> [1]
  y -> [1]
  z -> [1]

julia> R,(x,y) = polynomial_ring(base_ring(C8), [:x, :y])
(Multivariate polynomial ring in 2 variables over QQ, QQMPolyRingElem[x, y])

julia> C8tmp = C8((inv(B)*gens(Rh))...);

julia> U, t = polynomial_ring(QQ, :t);

julia> h = evaluate(C8tmp, [t, one(U), zero(U)]);

julia> degree(h), Hecke.n_real_roots(h)
(8, 0)

julia> zchart = hom(Rh, R, [x,y,1])
Ring homomorphism
  from graded multivariate polynomial ring in 3 variables over QQ
  to multivariate polynomial ring in 2 variables over QQ
defined by
  x -> x
  y -> y
  z -> 1

julia> draw_curve_tikz("shustin.tikz", zchart(C8tmp); solver_precision=1024, ybox_tolerance=512)

julia> draw_curve_tikz("shustin_topological.tikz", zchart(C8tmp); solver_precision=1024, ybox_tolerance=512, graph=true)

julia> QC4 = load("QC4.mrdi"; params=Rh);

julia> QC4tmp = QC4((inv(B)*gens(Rh))...);

julia> draw_curve_tikz("QC4.tikz", zchart(QC4tmp); solver_precision=1024, ybox_tolerance=512)
\end{minted}

The session first checks that the binary form $C_8(x,y,0)$ has no real zeros, so that $V_\RR(C_8)$ lies in the chosen affine chart.
The schematic rendering it produces is shown in Fig.~\ref{fig:shustin_topological}: four empty ovals and three ovals containing five empty ovals each.
These are $4+3\cdot 6=22$ ovals, which is the maximal number $M=\tfrac{1}{2}\cdot 7\cdot 6+1=22$ for $d=8$.
For this curve, the construction, including $Q(C_4)$ and its singularities, as well as the decisions of the drawing algorithm were checked by hand, and the topology was confirmed independently by Fabrice Rouillier with \pace.
This proves the theorem.

\section{Outlook}
The current implementation of plane curve visualization within \oscar serves as an entry point toward studying real plane algebraic curves.
As pointed out already, there is a multitude of papers containing various improvements that should and will be integrated into \oscar in the future.
On the software side, the current implementation is being kept as modular as possible, to allow for different backends, i.e., solvers, and other frontends, like \julia's \plotsjl framework.

\subsection*{Acknowledgments}
We would like to thank Elias Tsigaridas for a fruitful exchange at the MEGA conference and Fabrice Rouillier for independently verifying the topology of $C_8$ with \pace.
The authors received funding by the Deutsche Forschungsgemeinschaft (DFG, German Research Foundation) under Germany's Excellence Strategy -- \enquote{The Berlin Mathematics Research Center MATH$^+$} (EXC-2046/1, EXC-2046/2, pro\-ject ID 390685689), \enquote{Symbolic Tools in Mathematics and their Application} (TRR 195, project ID 286237555), and the \enquote{Mathematical Research Data Initiative (MaRDI)} (project ID 460135501).

\printbibliography
\end{document}